\documentclass{amsart}
\usepackage{graphicx}
\usepackage{mathtools}

\usepackage{amsmath, amsthm, amsfonts}
\usepackage{amssymb}
\usepackage{tikz}
\usetikzlibrary{matrix}

\newtheorem{thm}{Theorem}[section]

\newtheorem{prop}[thm]{Proposition}
\theoremstyle{definition}
\newtheorem{defn}[thm]{Definition}

\theoremstyle{remark}

\numberwithin{equation}{section}
\newcommand{\norm}[1]{\left\Vert#1\right\Vert}

\DeclareMathOperator*\uplim{\overline{lim}}
\begin{document}
\title[divergent Fourier series]{divergent Fourier series with respect to biorthonormal systems in function spaces near  $L^1$}%
\author{Nikoloz Devdariani}
\address[Nikoloz Devdariani]{Faculty of Exact and Natural Sciences, Javakhishvili Tbilisi State University, 13, University St., Tbilisi, 0143, Georgia}%
\email[Nikoloz Devdariani]{nikoloz.devdariani275@ens.tsu.edu.ge}%


\begin{abstract}
In this paper, we generalize Bochkarev's theorem, which states that for any uniformly bounded biorthonormal system $\Phi$, there exists a Lebesgue integrable function whose Fourier series with respect to the system $\Phi$ diverges on a set of positive measure. We find the class of variable exponent Lebesgue spaces $L^{p(\cdot)}([0,1]^n)$, where $1 < p(x) < \infty$ almost everywhere on $[0,1]^n$, such that the aforementioned Bochkarev's theorem holds.
\end{abstract}
\maketitle
\section{Introduction}
\par After Kolmogorov \cite{kolmogorov1}, \cite{kolmogorov2} presented examples of functions in $L^1$ with almost everywhere and everywhere divergent trigonometric Fourier series, many authors have attempted to generalize these results by providing examples of functions with almost everywhere divergent trigonometric Fourier series from narrower Orlicz spaces. The most significant result in this direction is due to Konyagin \cite{koniagin}, who achieved the same result for the space $L\varphi(L)$, provided that $\varphi$ satisfies $\varphi(t) = o(\sqrt{\ln t / \ln\ln t})$.
\par Similar problems with respect to other orthonormal systems have been considered by various authors. One such problem was posed by Alexits (see \cite[pp. 326]{alexits1}, \cite[pp. 287]{alexits2}) and Olevskii \cite{olevskii}, concerning an analogue of Kolmogorov's example of a divergent trigonometric Fourier series for general uniformly bounded orthonormal systems.
\par The answer to this question was provided by Bochkarev \cite{bochkarev}. He proved that for any given uniformly bounded orthonormal system, there exists a function in $L^1$ whose Fourier series with respect to this system diverges at every point of some set of positive measure. However, it turns out that a complete analogue of Kolmogorov's example for uniformly bounded orthonormal systems does not exist in general. This conclusion is based on Kazaryan's \cite{kazaryan} construction of a complete orthonormal system that is uniformly bounded and for which every Fourier series converges on some set of positive measure.
\par 
Later, Bochkarev \cite{bochkarev2} extended his aforementioned result to uniformly bounded biorthonormal systems defined on a separable metric space with a Borel regular outer measure.
\par 
The authors of the paper \cite{egk} provide a different perspective on the problem of almost everywhere divergence of trigonometric Fourier series in the subspaces of $L^1$, specifically in terms of variable Lebesgue spaces $L^{p(\cdot)}$. They show that $L^1 = \bigcup L^{p(\cdot)}$, where the union is taken over all measurable functions $p(\cdot)$ such that $p(x) > 1$ almost everywhere. This implies that any function whose Fourier series diverges almost everywhere must belong to some variable exponent space $L^{p(\cdot)}$ with $1 < p(x) < \infty$ almost everywhere.

In \cite{egk}, the authors construct a variable exponent space $L^{p(\cdot)}$, with $1 < p(x) < \infty$ almost everywhere, that shares with $L^\infty$ the property that the space of continuous functions $C$ is a closed linear subspace within it. Moreover, Kolmogorov's function, which has a Fourier series that diverges almost everywhere, belongs to the $L^{p'(\cdot)}$ space, where $p'(\cdot)$ is the Hölder conjugate of $p(\cdot)$. Additional results concerning the convergence of Fourier series for functions from these spaces can be found in \cite{k1} and \cite{k2}. Various other results related to these spaces are discussed in \cite{uribefiorenza} and \cite{DHHR}.
\par Later, in \cite{kopsamzv}, the authors provided an analogue of Bochkarev’s theorem for uniformly bounded orthonormal systems within a certain class of variable exponent Lebesgue spaces. They found the class of variable exponent Lebesgue spaces $L^{p(\cdot)}$, with $1 < p(x) < \infty$ almost everywhere, for which Bochkarev’s theorem holds.
\par For readers who wish to gain a deeper understanding of these subjects, we kindly recommend exploring the following works related to the research topic \cite{bochkarev*, Oniani}.
\par In the first of the aforementioned works, the almost everywhere divergence effect is established for any countable orthonormal system of characters of a compact group. In the second, the divergence effect is demonstrated everywhere for a broad class of character systems, including all Vilenkin systems.

\par Our plan for this paper is to characterize the class of variable exponent Lebesgue spaces for which an analogue of Bochkarev's theorem about bounded biorthonormal systems is valid. To achieve this, we will introduce some definitions and notations.
\begin{defn}
\label{definition of biorthogonal series}
Let $(X,S,\mu)$ be a measurable space, where $S$ is a $\sigma$-algebra of $\mu$-measurable sets and $\mu(X) = 1,$ and let $\{f_n,g_n\}$ be a biorthonormal system (see \cite[Ch. VIII, \S 1]{Kaczmarz}) such that $f_n,g_n \in L^{\infty}(X,\mu)$. For any function $F \in L^{1}(X,\mu)$ the system $\{f_n,g_n\}$  generates two Fourier series $$\sum_{n=1}^{\infty} (F,f_n)g_n(x) \: \: \: \: \text{and} \: \: \: \: \sum_{n=1}^{\infty} (F,g_n)f_n(x),$$
(These series are said to be conjugate).
\end{defn}
\par Bochkarev \cite{bochkarev2} proved the following
\begin{thm}
Let $X$ be a separable metric space with a Borel regular outer measure $\mu^*$ such
that $\mu^*(X) = 1$. Then, for any biorthonormal system $\{f_n, g_n\}$ satisfying conditions
\begin{equation}
\label{cond_syset_uniformly_boudedness}
\norm{f_n}_\infty,\norm{g_n}_\infty \leq A , \: \: n \in\mathbb{N},
\end{equation}
and
\begin{equation}
\label{cond_on_fourier_coefficients}
\lim_{n \to \infty} \int\limits_{E} f_n(x) d\mu(x) = 0 \: \: \text{and} \lim_{n \to \infty} \int\limits_{E} g_n(x)d\mu(x) = 0.
\end{equation}
for any measurable set $E \subset X$. Then there exist functions $F_1, F_2\in L^1(X, \mu)$ and a set $E\subset X$ such that 
$$
\mu (E) > 0
$$
and for all $x\in E$ we have
$$
\uplim_{N \to \infty} \sum_{n=1}^{N} ((F_1,g_n)f_n(x) + (F_2,f_n)g_n(x)) = \infty.
$$
\end{thm}
This theorem is valid, in particular, for the space $\mathbb{R}^n$, $n\in\mathbb{N}$ and any finite Borel regular
outer measure $\mu^*$. Thus, the following assertion holds (see \cite[Theorem 5]{bochkarev2})
\begin{thm}
\label{thm_bochkarev_particular_result}
If $\mu$ is the classical Lebesgue measure on $\mathbb{R}^{n}$ and $E$ is a measurable set with $\mu (E) < \infty$, then, for any uniformly bounded biorthonormal system $\{f_n,g_n\}$ on $E$ satisfying condition \eqref{cond_on_fourier_coefficients}, there exists a Fourier series divergent on a set of positive measure. 
\end{thm}
\par In this paper we are going to extend this Theorem \ref{thm_bochkarev_particular_result} for a variable Lebesgue spaces.
\par Let through the paper the symbol $|E|$ defines the Lebesgue measure of the measurable set $E$. Also, $\Omega := [0; 1]^{n}$ and for a given $p(\cdot)$, the Holder conjugate $p'(\cdot)$ is defined by $p'(x):=p(x)/(p(x)-1)$.
\par Let $(X,S,\mu)$ be a nontrivial measure space. Given a real-valued measurable function $f$ on $X$, define its decreasing rearrangement by
$$
f^*(s) = \inf\{\alpha \ge 0 \::\: \mu(\{x \in X : |f(x)| > \alpha\}) \le s\}, \quad s>0.
$$
Let $P_{\ln}$ be a set of all functions $p:X\to[1;\infty)$ such that
\begin{equation}
\label{cond_in_def_p_ln}
\limsup_{t\to0+}\frac{(p')^*(t)}{\ln(e/t)}>0.
\end{equation}
\par
Let $W(p)$ denote the set of all functions equimeasurable with $p(\cdot)$. Below, we will find the conditions on the function $p(\cdot)$ for which there exists $\bar{p}(\cdot) \in W(p)$ such that the space $C(\Omega)$ of continuous functions is a closed subspace in $L^{\bar p(\cdot)}(\Omega)$.
In this paper, we generalize the  result by Kopaliani, Samashvili and Zviadadze obtained in \cite{kopsamzv}, from the one-variable case to the several-variable case. Moreover, we extend this generalization to uniformly bounded biorthonormal systems. 
\par Let now state the result:
\begin{thm}
\label{thm_main}
For any biorthonormal system $\Phi:=\{f_n, g_n\}_{n\in\mathbb{N}}$ on $\Omega$, satisfying conditions \eqref{cond_syset_uniformly_boudedness} and \eqref{cond_on_fourier_coefficients}, and for any $p(\cdot) \in P_{\ln}$, there exists a measure preserving transformation $\omega : \Omega \to \Omega $, such that in the corresponding $L^{p(\omega(\cdot))}(\Omega)$ space    there exist functions $F_1, F_2 $ and a set $E \subset \Omega $ such that $|E| > 0$ and for all $x \in E$ we have
$$
\uplim_{N \to \infty} \sum_{n=1}^{N} ((F_1,g_n)f_n(x) + (F_2,f_n)g_n(x)) = \infty.
$$
\end{thm}
\section{definitions and auxiliary results}
\par Let $\mathcal{M}$ denote the space of all equivalence classes of Lebesgue measurable real-valued functions on $\Omega$, equipped with the topology of convergence in measure relative to each set of finite measure.
\begin{defn}
\label{defn_BFS}
A Banach subspace $X$ of $\mathcal{M}$ is referred to as a Banach function space (BFS) on $\Omega$ if the following conditions hold:

$1)$ The norm $\|f\|_{X}$ is defined for every measurable function $f$, and $f\in X$ if and only if $\|f\|_{X} < \infty$. Also, $\|f\|_{X} = 0$ if and only if $f = 0$ almost everywhere;

$2)$ $\||f|\|_{X} = \|f\|_{X}$ for all $f\in X$;

$3)$ If $0 \leq f \leq g$ almost everywhere, then $\|f\|_{X} \leq \|g\|_{X}$;

$4)$ If $0 \leq f_{n} \uparrow f$ almost everywhere, then $\|f_{n}\|_{X} \uparrow \|f\|_{X}$;

$5)$ If $E$ is a measurable subset of $\Omega$ with finite measure ($|E|<\infty$), then $\|\chi_{E}\|_{X} < \infty$, where $\chi_{E}$ is the characteristic function of $E$;

$6)$ For every measurable set $E$ with finite measure ($|E|<\infty$), there exists a constant $C_{E} < \infty$ such that $\int_{E}f(t)dt \leq C_{E}\|f\|_{X}$.
\end{defn}

\par Now, let's introduce various subspaces of a BFS $X$:
\begin{itemize}
\item A function $f$ in $X$ has an absolutely continuous norm in $X$ if $\|f\cdot\chi_{E_n}\|_{X} \to 0$ whenever ${E_n}$ is a sequence of measurable subsets of $\Omega$ such that $\chi_{E_n} \downarrow 0$ almost everywhere. The set of all such functions is denoted by $X_A$;

\item $X_B$ is the closure of the set of all bounded functions in $X$;

\item A function $f \in X$ has a continuous norm in $X$ if, for every $x \in \Omega$, $\lim_{\varepsilon \to 0+}\|f\chi_{B(x,\varepsilon)}\|_X = 0$, where $B(x,\varepsilon)$ is a ball centered at $x$ with radius $\varepsilon$. The set of all such functions is denoted by $X_C$.
\end{itemize}

\par The relationship between the concepts of $X_A$ and $X_B$ is given in \cite{BS}. Generally, the interplay among the subspaces $X_A$, $X_B$, and $X_C$ can be intricate. For instance, there exists a BFS $X$ in which ${0} = X_A \subsetneq X_C = X$ (for example see  \cite{LN}).

\par Let $\mathcal{P}$ through whole paper denotes the family of all measurable functions $p: \Omega\rightarrow[1;+\infty)$.
When $p(\cdot)\in\mathcal{P}$ we denote by $L^{p(\cdot)}(\Omega)$, the set of all measurable functions $f$ on $\Omega$ such that for some
$\lambda>0$
$$
\int_{\Omega}\left(\frac{|f(x)|}{\lambda}\right)^{p(x)}dx<\infty.
$$
This set becomes a BFS when equipped with the
norm
$$
\|f\|_{p(\cdot)}=\inf\left\{\lambda>0:\,\,
\int_{\Omega}\left(\frac{|f(x)|}{\lambda}\right)^{p(x)}dx\leq1\right\}.
$$

\par The variable exponent Lebesgue spaces $L^{p(\cdot)}(\Omega)$ and the corresponding variable exponent Sobolev spaces $W^{k,p(\cdot)}$ are of significant interest due to their applications in fluid dynamics, partial differential equations with non-standard growth conditions, calculus of variations, image processing, and more (refer to \cite{uribefiorenza,DHHR} for further details).

\par For the specific case of a particular BFS $X = L^{p(\cdot)}(\Omega)$, the relationship between this space and its subspaces, namely, $X_A$, $X_B$, and $X_C$, has been explored in \cite{ELN}. We will now present some of the key findings from that paper.
\par
\begin{prop}[Edmunds, Lang, Nekvinda]
\label{prop_X_A_equals_X_B}
Let $p(\cdot)\in \mathcal{P}$ and set $X=L^{p(\cdot)}(\Omega)$. Then 
\par (i) $X_A=X_C$;
\par (ii) $X_B=X$ if and only if $p(\cdot)\in L^\infty(\Omega)$;
\par (iii) $X_A=X_B$ if and only if
$$
\int_0^1 c^{p^*(t)}dt<\infty, \:\:\: \text{for all}\:\:\: c>1.
$$
\end{prop}
\par If $\psi$ is an increasing convex function $\psi:[0;+\infty)\to[0;+\infty)$, such that $\psi(0)=0$,
$$
\lim_{x\to0+}(\psi(x)/x)=0,\quad \text{and} \quad \lim_{x\to+\infty}(\psi(x)/x)=+\infty,
$$
then the Orlicz space $L_\psi$ is defined as the set of all $f\in \mathcal{M}(\Omega)$ for which:
$$
||f||_{L_\psi}=\inf\left\{\lambda>0\::\:\int_{\Omega}\psi\left(\frac{|f(t)|}{\lambda}\right)dt\leq1\right\}<+\infty.
$$

\par Recall that a nonnegative function $\varphi$ defined on $[0;+\infty)$ is called quasiconcave if it satisfies the following conditions: $\varphi(0)=0$, $\varphi(t)$ is increasing, and $\varphi(t)/t$ is decreasing.

\par The Marcinkiewicz space $M_\varphi$ is defined as the set of all $f\in \mathcal{M}(\Omega)$ for which:
$$
||f||_{M_\varphi}=\sup_{0<t}\frac{1}{\varphi(t)}\int_0^tf^*(u)du<+\infty.
$$
\par It is worth noting that $(M_\varphi)_A = (M_\varphi)_B$, and $(M_\varphi)_A$ can be characterized as the set of functions $f\in\mathcal{M}$ (see \cite{kps}) that satisfy:
\begin{equation}
\label{estim_limit_decreasing_rearrangement_over_phi}
\lim_{t\to0+}\frac{1}{\varphi(t)}\int_0^tf^*(u)du=0.
\end{equation}

\par Additionally, when $\psi(t)=e^t-1$ and $\varphi(t)=t\ln(e/t)$, the corresponding Orlicz and Marcinkiewicz spaces coincide (see \cite{BS}), and we denote them as $e^L$ and $M_{\ln}$. Furthermore, it can be observed that (see \cite[Corollary 3.4.28]{EE}):
\begin{equation}
\label{estim_norm_equivalences}
||f||_{e^L}\asymp||f||_{M_{\ln}}\asymp\sup_{0<t\leq1}\frac{f^*(t)}{\ln(e/t)}.
\end{equation}

The following result was initially established in \cite{egk} for the single-variable case, and our goal now is to extend it to the multi-variable scenario. Since the proof of this statement can be easily derived from the one provided in \cite{egk}, we will omit it here.
\par
\begin{thm}
\label{thm_necess_and_suff_C_closedc_in_X}
Let $X$ be a BFS on $\Omega$. The space $C(\Omega)$ of continuous functions is a closed linear subspace of $X$ if and only if there exists a positive constant $c$ such that for every rectangle $I \subset \Omega$, we have
\label{estim_C_is_closed_in_X}
$$c\leq ||\chi_{I}||_X.$$
\end{thm}
\par
\begin{thm}
\label{thm_necess_suff_on_decreasing_rearrangement}
For the existence of $\bar{p}(\cdot) \in W(p)$ for which $C(\Omega)$ is a closed subspace in $L^{\bar{p}(\cdot)}(\Omega)$, it is necessary and sufficient that
\begin{equation}
\label{cond_necess_suffic}
\limsup_{t \to 0+}\frac{p^*(t)}{\ln(e/t)}>0.
\end{equation}
\end{thm}
\par
The forthcoming proof closely follows the framework presented in \cite{kop-zv}. However, since we encounter some differences when extending the proof from the one-dimensional case to multiple dimensions, we have chosen to provide the complete proof for the sake of clarity.
\par Necessity. Since the space $C(\Omega)$ is closed in $L^{p(\cdot)}(\Omega)$, then by Theorem \ref{thm_necess_and_suff_C_closedc_in_X} there exists positive constant $d$ such that $d\leq ||\chi_{I}||_{p(\cdot)}$ for all rectangles $I$.
this implies $X_A\neq X_B$. Then by Proposition \ref{prop_X_A_equals_X_B} there exists $c>1$ such that 
\begin{equation}
\label{estim_int_c_pow_p_rearrangement_infty}
\int_0^1c^{p^*(t)}dt=+\infty.
\end{equation}
\par Consider two cases: 
\par Case 1) $p^*(\cdot)\in e^L$. Since (\ref{estim_int_c_pow_p_rearrangement_infty}) holds then function $p^*(\cdot)$ does not have absolute continuous norm that is $p^*(\cdot)\in e^L\backslash\left(e^L\right)_A$. Then by (\ref{estim_norm_equivalences}) we get that $p^*(\cdot)\in M_{\ln}\backslash(M_{\ln})_A$ and by (\ref{estim_limit_decreasing_rearrangement_over_phi}) it is obvious that
$$
\limsup_{t\to0+}\frac{1}{t\ln(e/t)}\cdot\int_0^tp^*(u)du>0,
$$
finally using ones more (\ref{estim_norm_equivalences}) from the last estimation we get (\ref{cond_necess_suffic}).

\par Case 2) $p^*(\cdot)\notin e^L$. Then by (\ref{estim_norm_equivalences})
$$
\sup_{0<t\leq1}\frac{p^*(t)}{\ln(e/t)}=+\infty,
$$
consequently (\ref{cond_necess_suffic}) holds.
The necessity part of the theorem proved.
\par Sufficiency. Let (\ref{cond_necess_suffic}) holds. For all $t\in[0;1]$ define function $h(t)=\min\{p^*(t), \ln(e/t)\}$. It is obvious that in this case holds
$$
\limsup_{t\to0+}\frac{h(t)}{\ln(e/t)}>0,
$$
then there exists a sequence $t_k\downarrow0$, such that
$$\frac{h(t_k)}{\ln(e/t_k)}\geq d,\quad k\in\mathbb{N},$$
for some positive number $d$. Now choose subsequence $(t_{k_n})$ such that $2t_{k_{n+1}}<t_{k_n}$, for all natural $n$. Since $t_k\downarrow0$, we can always choose such subsequence, so without loss of generality we can assume that sequence $(t_k)$ is already such.
\par Let given function $f$ defined by
$$
f(t)=d\cdot\ln(e/t_k), \:\:\: t\in(t_{k+1};t_k], \:\:\: k\in \mathbb{N} \quad\textnormal{and} \quad f(t)=1, \:\:\:t\in(t_1;1].
$$
It is clear that $h(t)\geq f(t)$ for all $t\in[0;1]$. Now choose positive number $c$ such that $c>e^{1/d}$ then we get 
\begin{equation}
\label{estim_c_pow_p_rearrangement_infty}
\int_0^1 c^{h(t)}dt=+\infty.
\end{equation}
Indeed,
$$
\int_0^1c^{h(t)}dt\geq\int_0^1c^{f(t)}dt>\int_{t_{k+1}}^{t_k}c^{d\cdot\ln(e/t_k)}dt=
$$
$$
=(t_k-t_{k+1})\cdot e^{d\cdot\ln c\cdot\ln(e/t_k)}>\frac{t_k}{2}\cdot\left(\frac{e}{t_k}\right)^{d\cdot\ln c}\to+\infty,\quad k\to+\infty.
$$
Choose decreasing sequence $\{a_k\}_{k\in\mathbb{N}}$, such that 
$$
\int_{a_{k+1}}^{a_k}c^{h(t)}dt=1.
$$
By (\ref{estim_c_pow_p_rearrangement_infty}) such sequence always can be chosen. Now let $\Delta_k=[a_{k+1};a_k]$, and $\{r_k\::\:k\in\mathbb{N}\}$ is a countable dense set in $[0;1]$. Define $b_k=-a_{k+1}+r_k$. Now let $A_k:=\Delta_k+b_k=[r_k;r_k+a_k-a_{k+1}]$. Let $g_k(t)=h(t) \cdot \chi_{\Delta_k}(t)$, $k\in\mathbb{N}$. Define functions $p_k(t)$ by the induction:
$$
p_1(t)=g_1(t-b_1)\chi_{[0;1]}(t),
$$
$$
p_k(t)=\left(p_{k-1}(t)(1-\chi_{\Delta_k}(t-b_k))+g_k(t-b_k)\right)\cdot\chi_{[0;1]}(t),\quad k>1.
$$

It is clear that $h(t)$ is decreasing and therefore $p_k(t)\leq p_{k+1}(t)$, for all $t\in[0;1]$ and all $k\in\mathbb{N}$. Also for all $k\in\mathbb{N}$ we have
\begin{equation}
\label{estim_int_p_k}
\int_0^1p_k(t)dt\leq\int_0^1h(t)dt\leq\int_0^1\ln(e/t)dt = 2.
\end{equation}
Now define $q(\cdot)$ function by
$$
q(t)=\lim_{k\to+\infty}p_k(t),\quad t\in[0;1].
$$
By (\ref{estim_int_p_k}) we get that the function $q(\cdot)$ is almost everywhere finite.
By the construction it is clear that $q^*(t)\leq h(t)\leq p^*(t)$. Now by the well known result (see \cite[Theorem 7.5]{BS}) there exists measure preserving transformation $\omega:[0;1]\to[0;1]$ such that $q(t)=q^*(\omega(t))$. Now define $\hat{p}(\cdot)$ by $\hat{p}(t) := p^{*}(\omega(t))$, $t\in[0;1]$. Since $q^*(t)\leq p^*(t)$ it is obvious that $q^*(\omega(t))\leq p^*(\omega(t))$, then for all $t\in(0;1)$ we get following inequality
\begin{equation}
\label{estim_bar_p_by_tilde_p}
q(t)\leq \hat{p}(t).
\end{equation}
Now, construct an exponential function $\bar{p}:\Omega \to [1,\infty)$, for which  the space of continuous functions will be a closed subspace inside its corresponding variable exponent Lebesgue space. For this, let's define measure-preserving mapping $\rho:\Omega \to [0;1]$, with the following rule: Suppose that, $x=(x_1,...,x_n)\in\Omega$ and for every $i\in\{1,...,n\}$ index, the representation of its corresponding coordinates be the following: $x_i=0.a_{i1}a_{i2}a_{i3}...$, then 
$$
\rho(x)=0.a_{11}a_{21}...a_{n1}a_{12}a_{22}...a_{n2}...
\:.$$ 
This mapping, which we have mentioned above, is well-known from literature. Thus, we can define the function $\bar{p}(x)=\hat{p}(\rho(x))$. In order to complete the proof, we should verify that the space of continuous function will be a closed subspace in its corresponding $L^{\bar{p}(\cdot)}(\Omega)$ space. For this purpose, we should show that there  exists a positive number $K$ such that for every rectangle $I \in \Omega $, we have $||\chi_{I}||_{\bar{p}}\ge K$. Consider any number $c > 1$, as in view of the fact that the set of binary rational numbers is dense everywhere in the set of all real numbers, for this reason we can find an $n$-dimensional binary rectangle $I^d$   for this given $n$-dimensional rectangle such that, $I^d \subset I $. Then because of the properties of function $\rho$, $c > 1$, and by \eqref{estim_bar_p_by_tilde_p}, we get 
$$
\int_Ic^{\bar{p}(x)}dx\ge \int_{I^d}c^{\bar{p}(x)}dx = \int_{(I^d)'}c^{\hat{p}(t)}dt\ge\int_{(I^d)'}c^{q(t)}dt,
$$ 
where $(I^d)'$ denotes one-dimensional binary interval taken out from $[0;1]$, for which $\rho(I^d)=(I^d)'$.  
By the construction of $q(\cdot)$ there exists number $k_0$ such that $A_{k_0}\subset (I^{d})'$. Then we get 
$$
\int_{(I^{d})'} c^{{q}(t)} dt \ge \int_{A_{k_0}} c^{{q}(t)} dt \ge \int_{A_{k_0}} c^{p_{k_0}(t)} dt = 
$$
$$ =\int_{A_{k_0}} c^{g_{k_0}(t-d_{k_0})} dt = \int_{r_{k_0}}^{r_{k_0} +t_{k_0} - t_{k_0 + 1}} c^{h(t-d_{k_0})\cdot\chi_{\Delta_{k_0}}(t-d_{k_0})} dt =
$$
$$
= \int_{t_{k_0 + 1}}^{t_{k_0}} c^{h(t)} dt \geq 1.
$$

Now by the definition of the norm in variable Lebesgue space and by the above estimations we get that for all $n$-dimensional rectangles $I \subset \Omega$ we have $||\chi_I||_{\bar{p}}>1/c$. By the Theorem \ref{thm_necess_and_suff_C_closedc_in_X} we get the proof of sufficiency of the Theorem \ref{thm_necess_suff_on_decreasing_rearrangement}.

\section{Proof of Theorem \ref{thm_main}}
\par Before proceeding with the construction of functions whose Fourier series diverge over sets of positive measure, we will establish foundational considerations. These considerations will serve as the basis upon which we will later develop the corresponding functions. It is important to note that we will not explicitly prove the divergence of the Fourier series at this point. Instead, our focus will be on constructing these functions, followed by outlining a method for applying Bochkarev's theorem to establish their divergence.

\par Let's start by constructing the $p(\cdot)$. For all $t\in(0;1)$ define the function $h(t):=\min\{(p')^*(t), \ln(e/t)\}$. By \eqref{cond_in_def_p_ln} it is obvious that
$$
\limsup_{t\to0+}\frac{h(t)}{\ln(e/t)}>0.
$$
Then there exists a sequence $t_k\downarrow0$, such that
$$\frac{h(t_k)}{\ln(e/t_k)}\geq a,\quad k\in\mathbb{N},$$
for some positive number $a$. 
\par It is obvious that we can choose subsequence $t_{k_n}$ such that $2t_{k_{n+1}}<t_{k_n}$. Let's choose a positive number $c$ such that $c>e^{1/a}$, then we get 
\begin{equation}
\label{estim_c_pow_p_rearrangement_infty}
\int_0^1 c^{h(t)}dt>\int_{t_{k_{n+1}}}^{t_{k_n}}c^{a\cdot\ln(e/t_{k_n})}dt=
\end{equation}
$$
=(t_{k_n}-t_{k_{n+1}})\cdot e^{a\cdot\ln c\cdot\ln(e/t_{k_n})}>\frac{t_{k_n}}{2}\cdot\left(\frac{e}{t_{k_n}}\right)^{a\cdot\ln c}\to+\infty,\quad n\to+\infty.
$$
According to (\ref{estim_c_pow_p_rearrangement_infty}) and the fact that $t_k\downarrow0$ we can choose the subsequence $(t_{k_{n_m}})$ from $(t_{k_n})$ such that 
$$
\int_{t_{k_{n_{m+1}}}}^{t_{k_{n_m}}}c^{h(t)}dt\geq1,\quad m\in\mathbb{N}.
$$
So without loss of generality we can assume that sequence $(t_k)$ is already such that
\begin{equation}
\label{properties_of_t_k}
1<a\ln(e/t_1), \quad 2t_{k+1}<t_k,\quad \int_{t_{k+1}}^{t_k}c^{h(t)}dt\geq1, \quad k\in\mathbb{N}.
\end{equation}

\par Let $\{l_k\::\:k\in\mathbb{N}\}$ be a fixed dense set on $(0;1)$ (below we will choose $l_k$ by using biorthonormal system $\Phi$). Let $r_k$, $k\in\mathbb{N}$ is the following numeration of the table
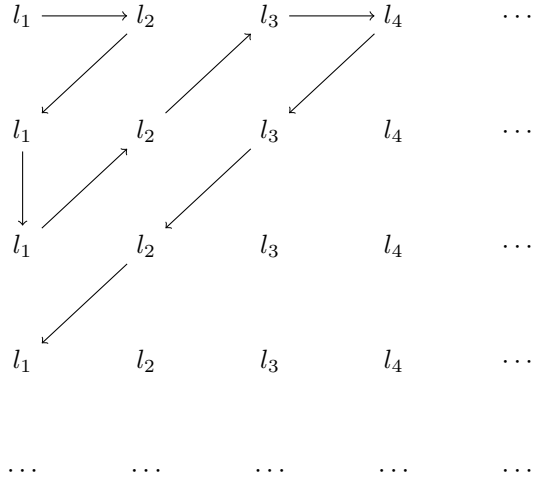
\begin{center}
\begin{tikzpicture}
\matrix(m)[matrix of math nodes,column sep=1cm,row sep=1cm]{
l_1 & l_2 & l_3 & l_4 & \cdots \\
l_1 & l_2 & l_3 & l_4 & \cdots \\
l_1 & l_2 & l_3 & l_4 & \cdots \\
l_1 & l_2 & l_3 & l_4 & \cdots \\
\cdots & \cdots & \cdots & \cdots & \cdots\\
};
\draw[->]
(m-1-1)edge(m-1-2)
(m-1-2)edge(m-2-1)
(m-2-1)edge(m-3-1)
(m-3-1)edge(m-2-2)
(m-2-2)edge(m-1-3)
(m-1-3)edge(m-1-4)
(m-1-4)edge(m-2-3)
(m-2-3)edge(m-3-2)
(m-3-2)edge(m-4-1);
\end{tikzpicture}
\end{center}
It is clear that for the each $l_k$ there exists sequence $(r_{k_m})$ such that $l_k=r_{k_m}$,  $m\in\mathbb{N}$.
Now let $\Delta_k:=[t_{k+1};t_k]$, where $t_k$ are points possessing the property (\ref{properties_of_t_k}). Define $d_k:=-t_{k+1}+r_k$ and $E_k:=\Delta_k+d_k=[r_k;r_k+t_k-t_{k+1}]$. Let $g_k(t):=h(t) \cdot \chi_{\Delta_k}(t)$, $k\in\mathbb{N}$. Let introduce the functions $q_k(t)$ by the induction:
$$
q_1(t):=g_1(t-b_1)\chi_{[0;1]}(t),
$$
$$
q_k(t):=\left[q_{k-1}(t)(1-\chi_{\Delta_k}(t-d_k))+g_k(t-d_k)\right]\cdot\chi_{[0;1]}(t),\quad k>1.
$$

It is clear that $h(t)$ is decreasing and therefore $q_k(t)\leq q_{k+1}(t)$, for all $t\in[0;1]$ and all $k\in\mathbb{N}$. Also for all $k\in\mathbb{N}$ we have
\begin{equation}
\label{estim_int_q_k}
\int_0^1q_k(t)dt\leq\int_0^1h(t)dt\leq\int_0^1\ln(e/t)dt = 2.
\end{equation}
Now define a function
$$
\hat{q}(t)=\lim_{k\to+\infty}q_k(t),\quad t\in[0;1].
$$
It is clear that 
\begin{equation}
\label{inequality_q_q_k}
\hat{q}(t)\geq q_k(t)\geq a\ln(e/t_k), \:\:\: t\in E_k,\:k\in\mathbb{N}.
\end{equation}
By (\ref{estim_int_q_k}) we get that the function $\hat{q}(\cdot)$ is a.e. finite.
According to the construction it is clear that $\hat{q}^*(t)\leq h(t)\leq (p')^*(t)$. It follows from the well known result (see \cite[Theorem 7.5]{BS}) that there exists measure preserving transformation $\zeta:[0;1]\to[0;1]$ such that $\hat{q}(t)=\hat{q}^*(\zeta(t))$. Now define $\tilde{q}(\cdot)$ by $\tilde{q}(t)=(p')^*(\zeta(t))$. Since $\hat{q}^*(t)\leq (p')^*(t)$ it is obvious that $\hat{q}^*(\zeta(t))\leq (p')^*(\zeta(t))$, then for all $t\in(0;1)$ we get following inequality
\begin{equation}
\label{estim_hat_q_by_tilde_q}
\hat{q}(t)\leq \tilde{q}(t).
\end{equation}
\par  Now, as in the proof of (\ref{thm_necess_suff_on_decreasing_rearrangement}), we can construct an exponential function $\bar{q}:\Omega \to [1,\infty)$, for which  the space of continuous functions will be a closed subspace inside its corresponding variable exponent Lebesgue space. Then,we get 
$$
\int_Ic^{\bar{q}(x)}dx\ge \int_{I^b}c^{\bar{q}(x)}dx = \int_{(I^b)'}c^{\tilde{q}(t)}dt\ge\int_{(I^b)'}c^{\hat{q}(t)}dt,
$$ 
\par where $(I^b)'$ denotes one-dimensional dyadic interval taken out from $[0;1]$, for which $\rho(I^b)=(I^b)'$.  
By the construction of $\hat{q}(\cdot)$ there exists number $k_0$ such that $E_{k_0}\subset (I^{b})'$. Then we get 
$$
\int_{(I^{d})'} c^{\hat{q}(t)} dt \ge \int_{E_{k_0}} c^{\hat{q}(t)} dt \ge \int_{E_{k_0}} c^{q_{k_0}(t)} dt = 
$$
$$ 
=\int_{E_{k_0}} c^{g_{k_0}(t-d_{k_0})} dt = \int_{r_{k_0}}^{r_{k_0} +t_{k_0} - t_{k_0 + 1}} c^{h(t-d_{k_0})\cdot\chi_{\Delta_{k_0}}(t-d_{k_0})} dt =
$$
$$
= \int_{t_{k_0 + 1}}^{t_{k_0}} c^{h(t)} dt \geq 1.
$$

Now by the definition of the norm in variable Lebesgue space and by the above estimations we get that for all $n$-dimensional rectangles $I \subset \Omega$ we have $||\chi_I||_{\bar{q}(\cdot)}>1/c$. By the  \cite[Theorem 4]{dev} we get that the space of continuous function is a closed subspace in $L^{\bar{q}}(\Omega)$.

\par Consider the function $\bar{p}(\cdot)$ which is Holder conjugate of $\bar{q}(\cdot)$. It is clear that $\bar{p}(\cdot)$ is equmeasurable to $p(\cdot)$. Since $\Omega$ is a finite nonatomic measure space and $\bar{p}$ is measurable, \cite[Theorem 7.5]{BS} guarantees that there exists a measure-preserving transformation $\omega_1 : \Omega \to [0;1]$ such that $\bar{p}(x) = {\bar{p}}^{*}(\omega_1(x))$ for almost every $x \in \Omega$. Similarly, applying the theorem to $p$, there exists a measure-preserving transformation $\omega_2 : \Omega \to [0;1]$ such that $p(x) = p^{*}(\omega_2(x))$ for almost every $x \in \Omega$. Since any measure-preserving transformation between non-atomic finite measure spaces is a measure space isomorphism, $\omega_2$ is also a bijective(up to null sets), therefore its inverse $\overset {-1}{\omega_2}$ is also measure preserving.  Now if we define $\omega = \overset {-1}{\omega_2} \circ \omega_1$, we get that $\bar{p}(x) = p(\omega(x))$ for almost every $x \in \Omega$.  Hence, there  exists a  measure preserving transformation $\omega:\Omega\to\Omega$ such that $\bar{p}(x)=p(\omega(x))$.
\par Let $M_k:=\rho^{-1}(E_k)$, $k\in\mathbb{N}$ and $C>\ln(e/t_1)\cdot(a\ln(e/t_1)-1)^{-1}$, then by (\ref{inequality_q_q_k}) and (\ref{estim_hat_q_by_tilde_q}) it is obvious that
$$1<\bar{p}(x)\leq 1+\frac{C}{\ln(e/t_k)},\quad x\in M_k.$$

Since $\rho$ is measure preserving $\rho(M_k)=E_k$ and $|E_k|\asymp t_k$, then by the last estimation we obtain
\begin{eqnarray}
\label{estim_chi_e_k_norm}
||\chi_{M_k}||_{\bar{p}(\cdot)}\asymp |E_k|=t_k-t_{k+1}\asymp t_k.
\end{eqnarray}
Indeed, using \cite[Theorem 2.45]{uribefiorenza} we obtain
$$
\frac{1}{2}||\chi_{M_k}||_1\leq||\chi_{M_k}||_{\bar{p}(\cdot)}\leq 2||\chi_{M_k}||_{1+\frac{C}{\ln(e/t_k)}}\asymp t_k.
$$
\par Finally, by (\ref{estim_chi_e_k_norm}) we obtain
\begin{equation}
\label{estim_sum_chi_e_k_norm}
\left\|\sum_{k=1}^\infty a_k \chi_{M_k}\right\|_{\bar{p}(\cdot)}\asymp\left\|\sum_{k=1}^\infty a_k \chi_{M_k}\right\|_1.
\end{equation}
\par Recall that for each $k$ there exists a sequence of natural numbers $(k_m)$, $m\in\mathbb{N}$ such that $l_k=r_{k_m}$,  $m\in\mathbb{N}$. Thus, we can rewrite (\ref{estim_sum_chi_e_k_norm}) in the following form
\begin{equation}
\label{estim_sum_sum_chi_e_k_norm}
\left\|\sum_{k=1}^\infty\sum_{m=1}^\infty a_{k_m}\chi_{M_{k_m}}\right\|_{\bar{p}(\cdot)} \asymp 
\left\|\sum_{k=1}^\infty\sum_{m=1}^\infty a_{k_m}\chi_{M_{k_m}}\right\|_{1}.
\end{equation}

\par Next, we should select $\{l_k\in\mathbb{N}\}$ set by using the system $\Phi$ in such a way that after definition corresponding functions $F_1$ and $F_2$ and using \eqref{estim_sum_sum_chi_e_k_norm} we will obtain the proof of a theorem. Our construction of the functions $F_1$ and $F_2$ is analogous to the one from the work of Bochkarev  \cite{bochkarev2}.

\par For all $\theta\in\Omega$ consider sequence of the binary cubes $Q_m(\theta)\subset\Omega$ such that $\theta\in Q_m(\theta)$, $\forall m\in\mathbb{N}$ and diameter tends to zero. By Lebesgue differentiation theorem for any $f \in L^{1}(\Omega)$ we have that 
$$
\lim_{m\to\infty} \frac{1}{|Q(\theta)|}\int_{Q(\theta)}f(x)dx=f(\theta),
$$
for almost every $\theta \in \Omega$.
\par Let $G^{2N}$ denote the set of points $\theta^{(2N)} = (\theta_1^{(2N)},...,\theta_{2N}^{(2N)}) \in \Omega^{2N}$ for which 
$$
\lim_{m \to \infty} \frac{1}{|Q_m(\theta_{2i-1}^{(2N)})|} \int\limits_{Q_m(\theta_{2i-1}^{(2N)})} g_n(x)dx = g_n(\theta_{2i-1}^{(2N)})
$$
and
$$
\lim_{m \to \infty} \frac{1}{|Q_m(\theta_{2i}^{(2N)})|} \int\limits_{Q_m(\theta_{2i}^{(2N)})} f_n(x)dx = f_n(\theta_{2i}^{(2N)}),
$$
for all $n \in\mathbb{N}$ and $i \in\{1,...,N\}$. It is clear that $|G^{2N}|=1$. Consider the following set $\Theta':=\{\theta_i^{(2N)}\::\: i\in\{1,...,2N\},\:N\in\mathbb{N}\}$.
If this set is not dense in $\Omega$, we examine a countable set $\Theta$ such that $\Theta'\subset\Theta$ and $\Theta$ is dense in $\Omega$. Let $l_k$, $k\in\mathbb{N}$ is some numeration of $\Theta$.

\par For all fixed $N\in\mathbb{N}$ and fixed $\theta^{(2N)}_i$, $i\in\{1,...,2N\}$ there exists the sequence $r_{i_k}^{(2N)}$, $k\in\mathbb{N}$ such that $\theta^{(2N)}_i=r_{i_k}^{(2N)}$, $k\in\mathbb{N}$. Then by $E_i^{(2N)}$ we define the binary interval such that $\rho(\theta_i^{(2N)})\in E_i^{(2N)}$ and $|E_i^{(2N)}|\le t_{i_1}-t_{i_1+1}$ and let $M_i^{(2N)}:=\rho^{-1}(E_i^{(2N)})$.

\par For a sequence $\{N_n\}$ of positive integers and a decreasing sequence $\{\varepsilon_n\}$ of positive numbers (which we will specify below), we set 
$$
F_1(x) = \sum_{n=1}^{\infty} \frac{\varepsilon_n}{N_n} \sum_{i = 1}^{N_n} \frac{\chi_{M_{2i-1}^{(2N)}}(x)}{|M_{2i-1}^{(2N)}|},\quad
F_2(x) = \sum_{n=1}^{\infty} \frac{\varepsilon_n}{N_n} \sum_{i = 1}^{N_n} \frac{\chi_{M_{2i}^{(2N)}}(x)}{|M_{2i}^{(2N)}|}.
$$
Now for the proof of the existence of a set $E \subset \Omega $ such that $|E| > 0$ and for all $x \in E$ we have
$$
\uplim_{N \to \infty} \sum_{n=1}^{N} ((F_1,g_n)f_n(x) + (F_2,f_n)g_n(x)) = \infty.
$$
We simply need to replicate Bochkarev's proof step by step without making any changes. Consequently, the detailed proof will not be provided in this context.
\par Finally, by the definitions of the functions $F_1$, $F_2$, and \eqref{estim_sum_sum_chi_e_k_norm}, we have $F_1, F_2 \in L^{\bar{p}(\cdot)}(\Omega)$. As mentioned earlier, $\bar{p}$ is equimeasurable to $p$, implying the existence of a measure-preserving transformation $\omega: \Omega \to \Omega$ such that $\bar{p}(x) = p(\omega(x))$. Therefore, we conclude that $F_1, F_2 \in L^{p(\omega(\cdot))}(\Omega)$. This completes the proof.

\section*{Acknowledgement}
I am deeply grateful to the referees for their careful reading of my paper and for their helpful comments and remarks. I also appreciate their suggestion of references \cite{bochkarev*, Oniani} and thank them for their valuable advice.

\end{document}